# A new data fitting method for stretched Gaussian noise: stretched least square method

Wei Xu[1,2], Yingjie Liang[1,*], Wen Chen[1,*]

[1]Institute of Soft Matter Mechanics, College of Mechanics and Materials, Hohai University, Nanjing, Jiangsu 211100, China

[2]Department of Earth, Atmospheric, and Planetary Sciences, Purdue University, West Lafayette, Indiana 47907, USA

**Corresponding author:** Yingjie Liang, Email: liangyj@hhu.edu.cn

Wen Chen, Email: chenwen@hhu.edu.cn

**Abstract:** Stretched Gaussian distribution is the fundamental solution of the Hausdorff derivative diffusion equation and its corresponding stretched Gaussian noise is a widely encountered non-Gaussian noise in science and engineering. The least square method is a standard regression approach to fit Gaussian noisy data, but has distinct limits for non-Gaussian noise. Based on the Hausdorff calculus, this study develops a stretched least square method to fit stretched Gaussian noise by using the Hausdorff fractal distance as the horizontal coordinate. To better compare with the least square method, different high levels of stretched Gaussian noise are added to real values. Numerical experiments show the stretched least square method is more accurate specific in fitting stretched Gaussian noise than the least square method.

## 1. Introduction

In recent decades, fitting noisy data especially non-Gaussian noise has attracted growing attention [1,2]. Stretched Gaussian noise is a kind of typical non-Gaussian noise, which has widely been encountered in system identification, data simulation, and integrated circuit testing [3-5], just to mention a few. Parameters estimation is the most common statistical inference problem that always arises in physics and practical engineering [6], and the estimation method used to obtain the parameters should obviously depend on the specific problem encountered in application. The least square method, the most frequently used in data fitting and estimation, is a standard regression approach [7]. And under certain conditions the least square method will be the optimum method if the data errors obey the Gaussian distribution [8]. However, the least square method cannot get the best fitting results and may be biased or grossly inaccurate when non-Gaussian noise becomes larger [9]. Thus, this study presents a stretched least square method (Stretched-LSM) to fit the stretched Gaussian noise data.

Many fitting correction methods have been proposed to process the non-Gaussian noise, such as pseudo linear methods [10], Filter [11], and Fourier transform [12]. However most of the existing methods have low stability and large fitting error. Among them the wavelet threshold shrinkage [13] and the wavelet proportional shrinkage [14] are the most widely used wavelet methods, which have gained a repaid development. But the selection of threshold and proportion is not an easy task. Although many existing methods can deal with non-Gaussian noise, there are very few fitting methods for the stretched Gaussian noisy data. It is found that the stretched Gaussian noise obeys the stretched Gaussian distribution which can supply a gap of Gaussian distribution that inaccurately and ineffectively describes some statistical phenomena [15]. At present, stretched Gaussian distribution has been widely adopted to anomalous diffusion, turbulence, image processing, digital watermarking, synthetic aperture radar and ultrasonic image [16-19]. Especially, when describe the anomalous diffusion process in fractal structural porous media, the superiorities of this statistical method can best embody [20, 21]. In fact, the stretched Gaussian distribution is the fundamental solution of the Hausdorff derivative diffusion equation. The Hausdorff derivative underlies to the Hausdorff fractal distance, which can be considered a scale transform of Euclidean distance [22].

In this paper, we employ the Hausdorff fractal distance to the horizontal coordinate and develop a Stretched-LSM to effectively fit the stretched Gaussian noisy data. The Stretched-LSM improves the least square method. It should be pointed out that some variations of the least square method transform the ordinates or cumulative errors into different forms [23]. While the main feature of the Stretched-LSM is to use the Hausdorff fractal distance as the horizontal coordinate. As we know, the observed data are often contaminated by various noises [24]. To examine the Stretched-LSM, the stretched Gaussian random numbers are generated as the noise, and then the noises are added to the exact values of the selected functions. The two selected functions are polynomial and trigonometric functions, which are used as the longitudinal coordinates of the observed values in experiment. Based on the scaling transform, the horizontal coordinate can be changed by varying the Hausdorff fractal dimensions. Then the observed values are respectively fitted by using the least square method and the Stretched-LSM, respectively. In order to check the fitting stability and accuracy, the maximum absolute and the mean square errors are calculated to compare the results of the Stretched-LSM fitting curves and the least square method fitting curves to those of the primitive functions.

The rest of the paper is organized as follows. Section 2 proposes the Hausdorff derivative methodologies to fit the stretched Gaussian noise data, and followed by several case studies in Section 3. The results are then discussed in Section 4. Finally, a brief summary is concluded upon the foregoing results and discussion.

**2. Methodologies**

*2.1 Hausdorff derivative and stretched Gaussian distribution*

Supposing that a particle moves uniformly along a curve in terms of fractal time, the movement distance $l$ can be defined as [22]:

$$l(\tau) = v(\tau - t_0)^{\alpha}, \tag{1}$$

where $v$ represents the uniform velocity, $\tau$ the current time instance, $t_0$ the initial instance, $\alpha$ the fractal dimensionality in time. For the variable speed motion, the Hausdorff integral distance is calculated by:

$$l(t) = \int_{t_0}^{t} v(\tau) \, d(\tau - t_0)^{\alpha}. \tag{2}$$

Then the Hausdorff derivative in time is given by:

$$\frac{dl}{dt^{\alpha}} = \lim_{t' \to t} \frac{l(t) - l(t')}{(t - t_0)^{\alpha} - (t' - t_0)^{\alpha}} = \frac{1}{\alpha (t - t_0)^{\alpha - 1}} \frac{dl}{dt}. \tag{3}$$

Let the initial instance $t_0 = 0$, Eq. (3) is reduced to the basic concept of Hausdorff derivative [15]:

$$\frac{dl}{dt^{\alpha}} = \lim_{t \to t'} \frac{l(t) - l(t')}{t^{\alpha} - t'^{\alpha}} = \frac{1}{\alpha t^{\alpha - 1}} \frac{dl}{dt}. \tag{4}$$

Similarly the Hausdorff derivative in space with the initial instance $x_0 = 0$ is stated as:

$$\frac{dl}{dx^{\beta}} = \lim_{x \to x'} \frac{l(x) - l(x')}{x^{\beta} - x'^{\beta}} = \frac{1}{\beta x^{\beta - 1}} \frac{dl}{dx}. \tag{5}$$

According to Eq. (4) and Eq. (5), the Hausdorff diffusion equation is proposed to characterize anomalous diffusion [15]:

$$\frac{\partial u}{\partial t^{\alpha}} = D\left(\frac{\partial}{\partial x^{\beta}}\left(\frac{\partial u}{\partial x^{\beta}}\right)\right). \tag{6}$$

Here $0 < \alpha \leq 1, 0 < \beta \leq 1$, $D$ represents diffusion coefficient. And the solution of Eq. (6) has the form of stretched Gaussian distribution, whose probability density function is given by:

$$u(x,t) = \frac{1}{\sqrt{4\pi D t^{\alpha}}} \exp\left(-\frac{x^{2\beta}}{4Dt^{\alpha}}\right), \tag{7}$$

when $\beta = 1$ and $\alpha = 1$, it reduces to the standard Gaussian distribution.

In fact, Eq. (6) can be deemed to a scale transform of the normal diffusion equation, by using the Hausdorff space-time distance, i.e., the fractal metric space-time:

$$\begin{cases} \Delta \hat{t} = \Delta t^{\alpha} \\ \Delta \hat{x} = \Delta x^{\beta} \end{cases}. \tag{8}$$

The method [25] to generate the stretched Gaussian distribution random numbers is mature. After normalization, they can be used as stretched Gaussian noise. In this study, we use the acceptance rejection method [26].

*2.2 The stretched least square method (Stretched-LSM)*

The problem of detecting constant signals in additive noise is described by the following hypothesis tests:

$$f : R^m \to R, \tag{9}$$

$$y_i = f(x_i) + \varepsilon_i \ (\varepsilon_i = noise \times \eta\%, \ i = 1, 2, ..., n), \tag{10}$$

where $n$ is the number of observations, $y_i \in R$ the response variable, $x_i \in R^m$ the explanatory variable, $\varepsilon$ means the noise that satisfies stretched Gaussian distribution and $\eta$ is the noise percentage. In this study, we consider the case $n = 200$, then the observed values are $y_1, y_2,..., y_{200}$. According to the Hausdorff space-time distance in Eq. (8), the procedure is to:

(1) Resetting the horizontal: $xx = x + x^\beta$, $\beta \leq 1$ is the fractal dimension.

(2) Calculating the transition function equation $F_T(x)$: using the least square method fitted the experimental value: $(xx_i, y_i)$.

(3) Establishing the regression equation $F(x)$: fitting transition point $(x_i, F_{Ti})$ through the least square method.

## 3. Several study cases and its results

To check the Stretched-LSM, various noise data are polluted to fit polynomial and trigonometric equations by adding different levels of stretched Gaussian noise to the exact values, finally the maximum absolute error and the mean square error are calculated for different cases and compared with the least square method. For convenience and clarification, some evaluation indicators are defined:

Error1: means Maximum absolute error: $\max |F(x_i) - f(x_i)|$.

Error2: mean square error: $\left\{ \sum_{i=1}^{n} \{F(x_i) - f(x_i)\}^2 / n \right\}^{0.5}$.

LSM: fitting results using the least square method.

Stretched-LSM: fitting results using the stretched least square method.

*3.1 Fitting polynomial curves by the stretched least square method*

A polynomial is constructed by means of addition, multiplication and exponentiation to a non-negative power, which is often used in engineering inspection and geological survey [27]. It can be written as the following form with a single variable $x$,

$$f(x) = a_n x^n + a_{n-1} x^{n-1} + ... + a_2 x^2 + a_1 x + a_0, \quad (a_n \neq 0), \tag{11}$$

where $a_0, a_1, ..., a_{n-1}, a_n$ are constants. Then we construct the following model:

$$f = x^2 + x + 2, \tag{12}$$

and the corresponding regression function,

$$F = ax^2 + bx + c, \quad (a \neq 0). \tag{13}$$

The fitting curves with different levels of noise and fractal dimension are depicted in Figs. 1-2 in the polynomial function cases with two different fractal dimension $\beta = 0.4, 0.8$. And Tables 1-4 give the corresponding estimated parameters and fitting errors. As shown in Figs. 1-2, the blue line fitting curves are closer than the red line fitting curves to the primitive functions with black line, which illustrates that the Stretched-LSM is superior to the least square method when the noise is larger in fitting noise data.

To well serve the method, in the polynomial cases, noise levels are set to 30% and 50%, and the fractal dimensions are 0.4 and 0.8. Tables 1-4 illustrate the coefficients for three different values and fitting errors. By applying the two methods to the stretched Gaussian noise data measured, it can be noticed that the maximum absolute and the mean square errors of Stretched-LSM are smaller than those of the least square method, with the relationship expressed as:

$$\begin{cases} \text{Error1}(Stretched\text{-}LSM) < \text{Error1}(LSM) \\ \text{Error2}(Stretched\text{-}LSM) < \text{Error2}(LSM) \end{cases}. \tag{14}$$

By using the Stretched-LSM, the range of fitting maximum absolute error is (0.0028, 0.0194), the minimum of mean absolute error is 0.0021, and 0.075 is the largest. Meanwhile there is a regular, the bigger the noise, the worse the fitting. But for the least square method, maximum absolute error scope is (0.0089, 0.0309), the minimum of mean absolute error is 0.0051, and 0.0118 is the biggest. In the experiments, there are some special cases that the mean square errors of the two fitting methods are same, but the absolute error of the least square method is large, and that of Stretched-LSM is relatively small.

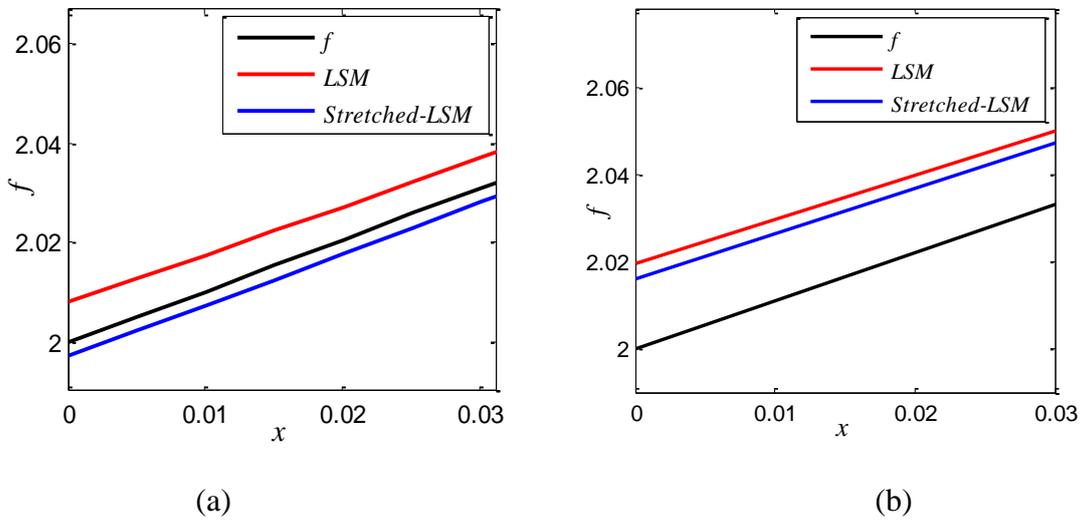

Fig. 1. Polynomial model fitting to the two cases $\eta = 30$: (**a**) $\beta = 0.4$; (**b**) $\beta = 0.8$.

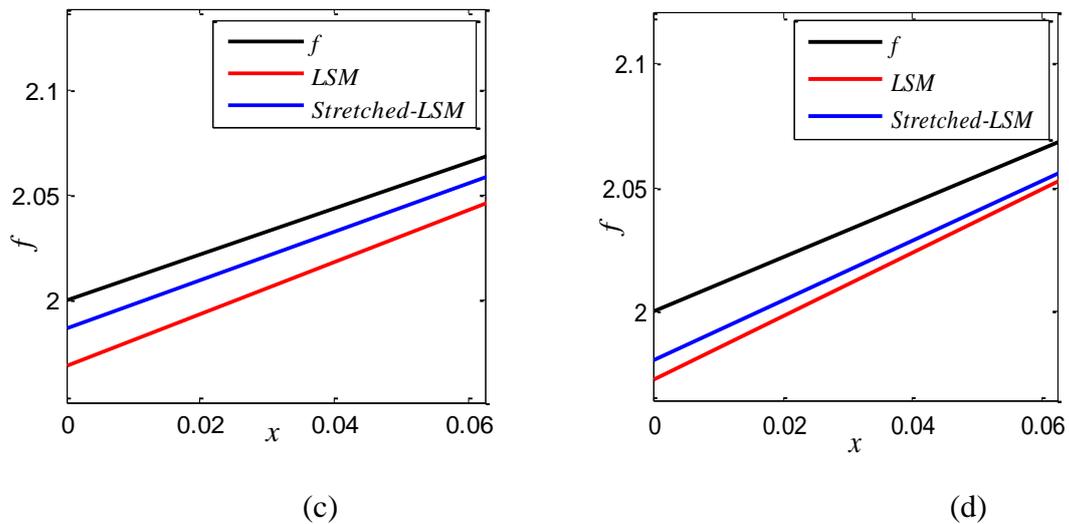

Fig. 2. Polynomial model fitting to the two cases $\eta = 50$: (**c**) $\beta = 0.4$; (**d**) $\beta = 0.8$.

Table 1. The estimated results for $\beta = 0.4$ and $\eta = 30$ in the polynomial case.

| Title 1 | a | b | c | Error1 | Error2 |
|---|---|---|---|---|---|
| f | 1 | 1 | 2 | 0 | 0 |
| LSM | 1.0625 | 0.9385 | 2.0079 | 0.0089 | 0.0051 |
| Stretched-LSM | 0.5497 | -0.1176 | 2.0275 | 0.0028 | 0.0021 |

Table 2. The estimated results for $\beta = 0.8$ and $\eta = 30$ in the polynomial case.

| Title 2 | a | b | c | Error1 | Error2 |
|---|---|---|---|---|---|
| f | 1 | 1 | 2 | 0 | 0 |
| LSM | 1.0770 | 0.9109 | 2.0194 | 0.0194 | 0.0068 |
| Stretched-LSM | 0.3388 | 0.3148 | 2.0165 | 0.0162 | 0.0056 |

Table 3. The estimated results for $\beta = 0.4$ and $\eta = 50$ in the polynomial case.

| Title 3 | a | b | c | Error1 | Error2 |
|---|---|---|---|---|---|
| f | 1 | 1 | 2 | 0 | 0 |
| LSM | 0.8734 | 1.1497 | 1.9691 | 0.0309 | 0.0118 |
| Stretched-LSM | 0.5432 | -0.0921 | 2.0112 | 0.0128 | 0.0065 |

Table 4. The estimated results for $\beta = 0.8$ and $\eta = 50$ in the polynomial case.

| Title 4 | a | b | c | Error1 | Error2 |
|---|---|---|---|---|---|
| f | 1 | 1 | 2 | 0 | 0 |
| LSM | 0.6098 | 1.2172 | 1.9726 | 0.0274 | 0.0099 |
| Stretched-LSM | 0.2706 | 0.4230 | 1.9746 | 0.0194 | 0.0075 |

*3.2 Fitting Trigonometric function by the stretched least square method*

Trigonometric function is widely applied to signal transmission, detection and estimation [28]. In recent years, trigonometric function as a new blind signal separation algorithm is proposed to get rid of the influence of surrounding noise [29]. In this part, we select a simple function:

$$f(x) = \sin(x), \tag{15}$$

The hypothetical regression equation is:

$$F = a\sin(bx+c)+d, (a \neq 0). \tag{16}$$

The same as polynomial fitting, in the trigonometric function case, noise levels are set to 30% and 50%, and the fractal dimensions are 0.4 and 0.8. Figs. 3-6 illustrate that the fitting curves with different levels of noise and fractal dimension. It shows that the blue line fitting curves are closer than the red line fitting curves to the primitive functions with pink line.

Tables 5-8 give the corresponding estimated parameters for four different values and fitting errors. Compared with the errors, by using the Stretched-LSM the range of fitting maximum absolute error is (0.0211, 0.0377), and the minimum of mean absolute error is 0.0066, and 0.0126 is the biggest. But for the least square method, maximum absolute error scope is (0.0411, 0.0633), the minimum of mean absolute error is 0.0159, and 0.0187 is the biggest.

By applying the two methods to the stretched Gaussian noise data measured, the maximum absolute and the mean square errors of the Stretched-LSM are smaller than those of the least square method, the relationship is also subordinated to Eq. (14). And from Figs. 1 to 6, it can be summarized, that the Stretched-LSM fitting curves are closer to the primitive functions than the least square method fitting curves. Thus, for different fractal dimensions, noise levels and different equation types, the Stretched-LSM is more accurate when it is applied to the stretched Gaussian noise data fitting compared with the least square method, especially the noise level is larger.

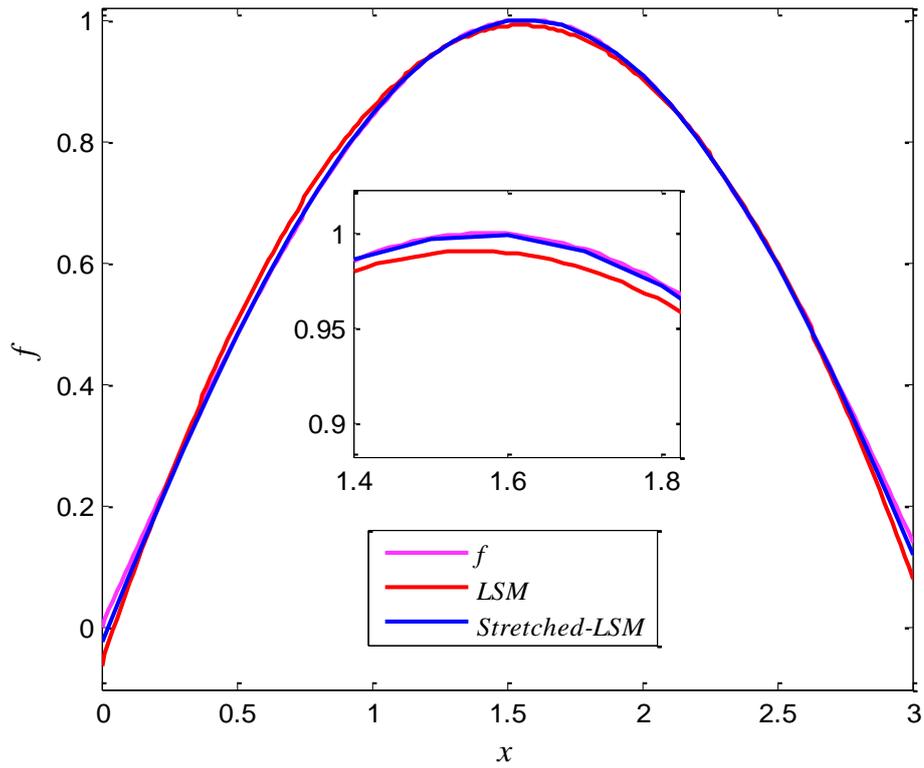

Fig. 3. Trigonometric function model fitting to the two cases $\eta = 30$, $\beta = 0.4$.

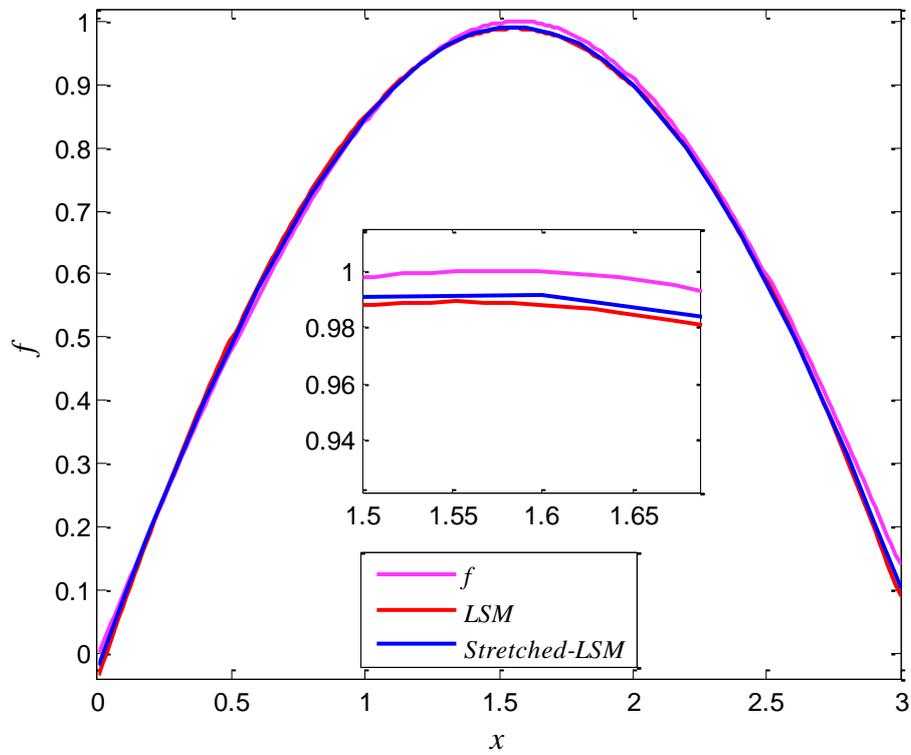

Fig. 4. Trigonometric function model fitting to the two cases $\eta = 30$, $\beta = 0.8$.

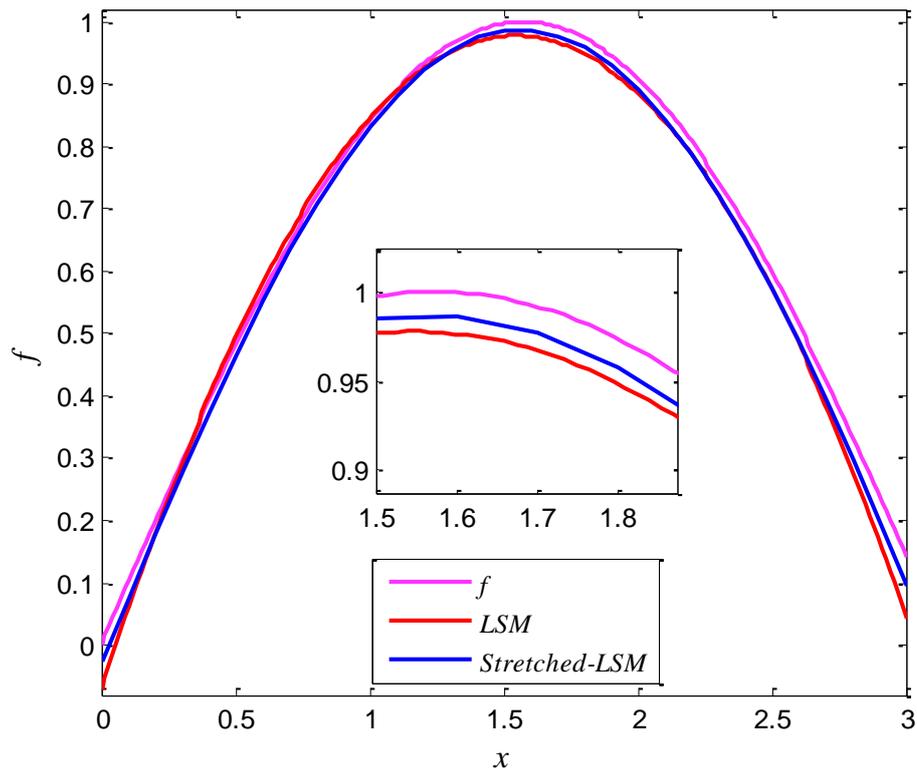

Fig. 5. Trigonometric function model fitting to the two cases $\eta = 50$, $\beta = 0.4$.

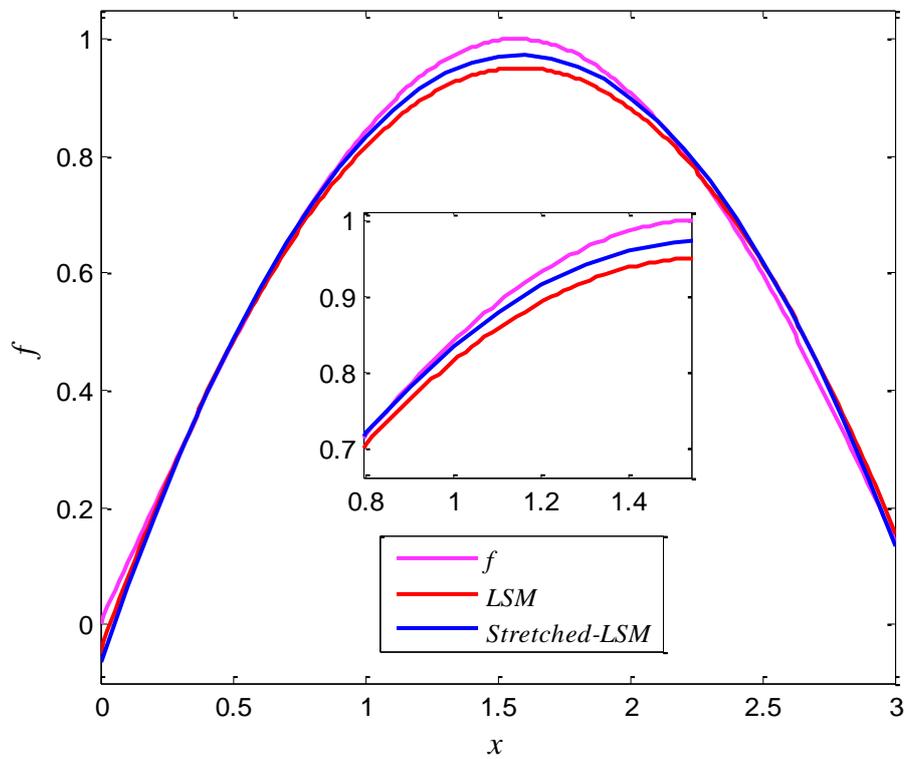

Fig. 6. Trigonometric function model fitting to the two cases $\eta = 50$, $\beta = 0.8$.

Table 5. The estimated results for $\beta = 0.4$ and $\eta = 30$ in the trigonometric function case.

| Title 5 | a | b | c | d | Error1 | Error2 |
|---|---|---|---|---|---|---|
| f | 1 | 1 | 0 | 0 | 0 | 0 |
| LSM | 18.4621 | 0.2182 | 1.2313 | -17.4714 | 0.0633 | 0.0187 |
| Stretched-LSM | 1.2292 | 0.8941 | 0.1708 | -0.2297 | 0.0211 | 0.0066 |

Table 6. The estimated results for $\beta = 0.8$ and $\eta = 30$ in the trigonometric function case.

| Title 6 | a | b | c | d | Error1 | Error2 |
|---|---|---|---|---|---|---|
| f | 1 | 1 | 0 | 0 | 0 | 0 |
| LSM | 2.2898 | 0.6344 | 0.5855 | -1.3009 | 0.0510 | 0.0161 |
| Stretched-LSM | 1.4188 | 0.8227 | 0.2907 | -0.4265 | 0.0377 | 0.0126 |

Table 7. The estimated results for $\beta = 0.4$ and $\eta = 50$ in the trigonometric function case.

| Title 7 | a | b | c | d | Error1 | Error2 |
|---|---|---|---|---|---|---|
| f | 1 | 1 | 0 | 0 | 0 | 0 |
| LSM | 3.9744 | 0.4745 | 0.8331 | -2.9825 | 0.0480 | 0.0159 |
| Stretched-LSM | 0.7701 | 1.1773 | -0.2781 | 0.2368 | 0.0253 | 0.0068 |

Table 8. The estimated results for $\beta = 0.8$ and $\eta = 50$ in the trigonometric function case.

| Title 8 | a | b | c | d | Error1 | Error2 |
|---|---|---|---|---|---|---|
| f | 1 | 1 | 0 | 0 | 0 | 0 |
| LSM | 2.6108 | 0.5915 | 0.6521 | -1.6119 | 0.0411 | 0.0161 |
| Stretched-LSM | 1.4618 | 0.8072 | 0.3125 | -0.4603 | 0.0221 | 0.0107 |

**4. Discussion**

The main ideas of this study are to change the horizontal coordinate using Hausdorff fractal space distance, then to obtain the transition equation and transition points by using the least square method, finally to get the fitting result based on the least square method. In particular, we generate random numbers with good randomness and obeying stretched Gaussian distribution instead of noise. From the above content, the cases with two parameters were detected: the noise percentage $\eta$ and the stretched exponent $\beta$. For the number of observations $n$, there is no effect on the fitting results at the tests. This indicates that the number of noise and the fitting results are not correlated. Meanwhile, the stretched Gaussian noise is closely related to the fractal derivative which was defined on non-Euclidean fractal metrics by using a time-space scaling transform [30]. Fractal theory has been applied in many aspects at present [31-33], this paper employs it on noise processing not only provides reliable theory bases and technology preparation for the application but also enriches its studying and applying field.

Because of the objective reasons such as environment, time and climate, the original data collected in many industrial situations are mixed with a large number of noise data [34]. In this paper, the test data are processed by curve fitting based on polynomial function and trigonometric function with periodicity which are fundamental models for several biomedical applications. To increase the confidence of the stretched least square method, we randomly selected 100 sets noise data for each test under identical conditions. The results from fitting error calculated that about 70% results are consistent with Eq. (14) and 30% of the remaining data displays the least square method fits well. In the real environmental engineering, noises are generated by random vibration, even if in the same distribution the random numbers are different, some are peak value, and some are empennage. There are number of deviations but all of them are reasonable in the tests.

In order to reduce noise and judge the data accurately, many methods and techniques have been proposed [12, 15, 35]. At present, curve fitting technology [35] has an increasing trend and has great practical significance to obtain an approximate function expression between the measured physical points and the measured physical quantities. Moreover, different methods were applied in different occasions according to their own characters. The proposed method in this paper has important theoretical significance for the signal processing in complex and fractal systems. It should be pointed out that in this study the proposed method is only from the view of probability statistics without rigorous mathematical proof. For more complicated systems, it is necessary to prove the method mathematically which will be the focus of our future work. On the other hand, the precondition of using this method is to estimate the parameters of the noise, which is also a hot issue in non-Gaussian noise study.

## 5. Conclusion

In this study, the stretched least square method, a new fitting method, is proposed by changing the horizontal coordinate with Hausdorff fractal space distance. The method is tested by adding different high levels of stretched Gaussian noise to real values of the polynomial and exponential equations. The maximum absolute and the mean square errors of the stretched least square method and the least square method are calculated for the different cases. Based on the foregoing results and discussion, the following conclusions are drawn:

(1) The fitting results with the stretched least square method are more accurate than the least square method when stretched Gaussian noise level is large.

(2) The stretched least square method is an alternative mathematical method to fit the stretched Gaussian noise data.

(3) The stretched least square method has clear physical mechanism in the context of Hausdorff metrics.

**Acknowledgments**

The work described in this paper was supported by the National Natural Science Foundation of China (Nos. 11702085, 11772121), the Fundamental Research Funds for the Central Universities (No. 2019B16114), the China Postdoctoral Science Foundation (No. 2018M630500), and the China Scholarship Council (CSC) (Grant No. 201806710084).

**References**

[1] D. Jakubisin, R. Buehrer, Approximate joint MAP detection of Co-Channel signals in non-Gaussian noise, Trans. Commun. 99 (2016) 1-6.

[2] A. Ichiki, Y. Tadokoro, Relation between optimal nonlinearity and non-Gaussian noise: enhancing a weak signal in a nonlinear system, Phys. Rev. E 87 (1) (2013) 012124.

[3] F. Pascal, L. Bombrun, J.Y. Tourneret, Y. Berthoumieu, Parameter estimation for multivariate generalized Gaussian distributions, Trans. Signal Process. 61 (23) (2013) 5960-5971.

[4] A. Ichiki, Y. Tadokoro, Relation between optimal nonlinearity and non-Gaussian

noise: enhancing a weak signal in a nonlinear system, Phys. Rev. E 87 (1) (2013) 012124.

[5] E. Saatci, A. Akan, Respiratory parameter estimation in non-invasive ventilation based on generalized Gaussian noise models, Signal Process. 90 (2) (2010) 480-489.

[6] X.L. Li, J. G. Li, Y. M. Wang, The influence of non-Gaussian noise on the accuracy of parameter estimation, Phys. Lett. A 381 (4) (2016) 1-5.

[7] S.K. Oh, W. Pedrycz, A new approach to self-organizing fuzzy polynomial neural networks guided by genetic optimization, Phys. Lett. A 345 (1) (2005) 88-100.

[8] M. Karimzadeh, A.M. Rabiei, A. Olfat, Soft-Limited Polarity-Coincidence-Array spectrum sensing in the presence of non-Gaussian noise, IEEE Trans. Veh. Technol. 66 (2) (2017) 1418-1427.

[9] W. Xu, W. Chen, Y.J. Liang, Feasibility study on the least square method for fitting non-Gaussian noise data, Physica A 492 (2018) 1917-1930.

[10] F. Wang, F. Ding, Gradient-based iterative identification methods for multivariate pseudo-linear moving average systems using the data filtering, Nonlinear Dyn. 84 (4) (2016) 2003-2015.

[11] C. Bai, Y. Yang, K.W. Wei, A wavevector filter in inversion symmetric Weyl semimetal, Phys. Lett. A 380 (5-6) (2016) 764-767.

[12] J.W. Fan, Denoise in the pseudopolar grid Fourier space using exact inverse pseudopolar Fourier transform, Physics 5 (2015) 218-221.

[13] A. Antoniadis, D. Leporini, J. Pesquet, Wavelet thresholding for some classes of

non-Gaussian noise, Statist. Neerlandica 56 (4) (2002) 434-453.

[14] S. Didas, J. Weickert, Integrodifferential equations for continuous multiscale wavelet shrinkage, Inverse Probl. & Imag. 1 (1) (2017) 47-62.

[15] W. Chen, Time-space fabric underlying anomalous diffusion, Chaos Solitons Fract. 28 (4) (2005) 923-929.

[16] J.F. Bercher, Escort entropies and divergences and related canonical distribution, Phys. Lett. A 375 (33) (2011) 2969-2973.

[17] W. Chen, H. Sun, X. Li, The fractional derivative model of mechanics and engineering problems, Science Press, Beijing, 2010.

[18] Y.F. Zhu, A robust blind image watermarking based on generalized Gaussian distribution, Inform. Technol. J. 13 (7) (2014) 1427-1430.

[19] P. Saha, S. Banerjee, A.R. Chowdhury, Chaos, signal communication and parameter estimation, Phys. Lett. A 326 (1) (2004) 133-139.

[20] D.H. Zanette, P.A. Alemany, Thermodynamics of anomalous diffusion, Phys. Rev. Lett. 75 (1995) 366-369.

[21] W. Chen, Wang F, B. Zheng, W. Cai, Non-Euclidean distance fundamental solution of Hausdorff derivative partial differential equations, Eng. Anal. Bound. Elem. 84 (2017) 213-219.

[22] W. Chen, Fractal analysis of Hausdorff calculus and fractional calculus models, Compu. Aided Eng. 26 (3) (2017) 1-5.

[23] S.A. Hosseini, Neutron spectrum unfolding using artificial neural network and modified least square method, Radiat. Phys. Chem. 126 (2016) 75-84.

[24] H. Xian, Study on ANN noise adaptability in application of industry process characteristics mining, Adv. Mater. Res. 462 (2012) 635–640.

[25] D. Fulger, E. Scalas, G. Germano, Random numbers from the tails of probability distributions using the transformation method, Fract. Calc. Appl. Anal. 16 (2) (2013) 332-353.

[26] G. Horváth, M. Telek, Acceptance-Rejection Methods for generating random variates from matrix exponential distributions and rational arrival processes, ACM, New York, 2012.

[27] J. Kaklamanos, L.G. Baise, D.M. Boore, Estimating unknown input parameters when implementing the NGA Ground-Motion prediction equations in engineering practice, Earthq. Spectra 27 (4) (2011) 1219-1235.

[28] M. Syam, H.M. Jaradat, M. Alquran, A study on the two-mode coupled modified Korteweg-de Vries using the simplified bilinear and the trigonometric-function methods, Nonlinear Dyn. 90 (2) (2017) 1363-1371.

[29] S. Zhang, New exact solutions of the KdV-Burgers-Kuramoto equation, Phys. Lett. A 358 (5) (2017) 414-420.

[30] Y. Liang, W. Chen. A non-local structural derivative model for characterization of ultraslow diffusion in dense colloids. Commun Nonlinear Sci Numer Simulat 56 (2017) 131-137.

[31] Q. Wang, X. Shi, J. He, Z. Li. Fractal calculus and its application to explanation of biomechanism of polar bear hairs. Fractals 26 (2018) 1850086

[32] X. J. Yang, J. A. T. Machado, D. Baleanu. On exact traveling-wave solution for

local fractional Boussinesq equation in fractal domain. Fractals 25 (4) (2017) 1740006.

[33] Y. Wang, Q. Deng. Fractal derivative model for tsunami travelling. Fractals, 2018.

[34] P. R. PATNAIK. Fractal analysis of Lipase-Catalysed synthesis of butyl butyrate in a microbioreactor under the influence of noise. Fractals 21 (01) (2013) 1350007.

[35] I. I. Shevchenko, Lyapunov exponents in resonance multiplets, Phys. Lett. A 378 (1-2) (2014) 34-42.